\documentclass[11pt,leqno]{article}
\usepackage{amsmath, amssymb, amsthm}

\theoremstyle{plain}
\newtheorem{thm}{Theorem}
\newtheorem{cor}{Corollary}
\newtheorem{lem}{Lemma}[section]
\newcommand{\bvec}[1]{\mbox{\boldmath$#1$}}

\numberwithin{equation}{section}

\title{On the fourth power moment of $\Delta(x)$ and $E(x)$ \\
in short intervals}
\author{Yoshio Tanigawa  and  Wenguang Zhai}
\date{}

\begin{document}

\maketitle

\begin{abstract}
Let $\Delta(x)$ and $E(x)$ be error terms of the sum of divisor
function and the mean square of the Riemann zeta function,
respectively. In this paper their fourth power moments for short
intervals of Jutila's type are considered. We  get an asymptotic
formula for $U$ in some range.
\end{abstract}

 \footnote[0]{The second-named author is
supported by National Natural Science Foundation of China (Grant No.
10771127) and National Natural Science Foundation of Shandong
Province(Grant No. Y2006A31).}

%************************* section 1 *******************************************

\section{Introduction}

Let $d(n)$ denote the Dirichlet divisor function and $\Delta(x)$ denote the error term of
the sum $\sum_{n\leq x}d(n)$ for a large real variable $x$:
$$
\Delta(x)=\sum_{n \leq x}d(n)-x(\log x+2\gamma-1),
$$
where $\gamma$ is Euler's constant. Dirichlet first proved that
$\Delta(x)=O(x^{1/2})$. The exponent $1/2$ was improved by many authors.
The latest result reads
\begin{equation}
\Delta(x)\ll x^{131/416}(\log x)^{26947/8320 },
\end{equation}
proved by Huxley \cite{Hu2}. It is conjectured that
\begin{equation}
\Delta(x)=O(x^{1/4+\varepsilon}),
\end{equation}
which is supported by the classical mean square result
\begin{equation}
\int_1^T\Delta^2(x)dx=\frac{\zeta^4(3/2)}{6\pi^2\zeta(3)}T^{3/2} +O(T\log^5 T)
\end{equation}
proved by Tong \cite{To}. For results of higher power moments of $\Delta(x),$
see \cite{HB, I, IS, Ts, Z1, Z2, Z3, Z4}.

Define the function  $E(t)$  by
\begin{equation}
E(t):=\int_0^t\left|\zeta(\frac{1}{2}+iu)\right|^2du-t\log(t/2\pi)-(2\gamma-1)t,\hspace{2mm }t\geq 2.
\end{equation}
Many properties of $E(t)$ are parallel to those of $\Delta(x).$
Huxley \cite{Hu1} proved
\begin{equation} E(t)=O(t^{72/227}\log^{629/227}t).
\end{equation}
The conjectured bound is
\begin{equation}
E(t)=O(t^{1/4+\varepsilon}),
\end{equation}
which is supported by the mean square formula
\begin{equation}
\int_2^TE^2(t)dt=\frac{2\zeta^4(3/2)}{3\zeta(3)\sqrt{2\pi}}T^{3/2}
+O(T\log^5 T)
\end{equation}
proved by Meurman \cite{M}. For higher power moments of $E(t)$, see
\cite{HB, I, IS, Ts, Z1, Z2, Z3, Z4, Z5}.

Jutila \cite{J} first studied the mean square of the difference
$\Delta(x+U)-\Delta(x)$ for short intervals. He proved that if
$T\geq 2$ and $1\leq U\ll T^{1/2}\ll H\leq T,$ then
\begin{align} \label{eq18}
&\int_T^{T+H}\left(\Delta(x+U)-\Delta(x)\right)^2dx\\
& \qquad =\frac{1}{4\pi^2}\sum_{n\leq \frac{T}{2U}}\frac{d^2(n)}{n^{3/2}}\int_T^{T+H}x^{1/2}
   \left|\exp\left(2\pi i (n/x)^{1/2}U\right)-1\right|^2dx \nonumber \\
& \qquad +O(T^{1+\varepsilon}+HU^{1/2}T^{\varepsilon}),
\nonumber
\end{align}
which implies that the estimate
\begin{equation} \label{eq19}
\int_T^{T+H}\left(\Delta(x+U)-\Delta(x)\right)^2dx\asymp HU\log^3
(T^{1/2}/U)
\end{equation}
holds for $HU\gg T^{1+\varepsilon}$ and $T^{\varepsilon}\ll U\leq \frac{T^{1/2}}{2}.$

For $E(t)$ Jutila proved  that the asymptotic formula
\begin{align}  \label{eq110}
&\int_T^{T+H}\left(E(t+U)-E(t)\right)^2dt \\
&=\frac{1}{\sqrt{2\pi}}\sum_{n\leq T/2U}\frac{d^2(n)}{n^{3/2}}
\int_T^{T+H}t^{1/2}|\exp(i(2\pi n/t)^{1/2}U)-1|^2dt \nonumber \\[1ex]
& \quad +O(T^{1+\varepsilon}+HU^{1/2}T^{\varepsilon}), \nonumber
\end{align}
holds for  $1\leq U\ll T^{1/2}\ll H\leq T.$

Jutila \cite{J} also raised the problem of extending \eqref{eq18} and \eqref{eq110}
to higher power moments. Especially, he conjectured that the estimates
\begin{equation} \label{Jutila-conj}
\int_2^{T}\left(\Delta(x+U)-\Delta(x)\right)^4dx \ll
T^{1+\varepsilon}U^2
\end{equation}
and
\begin{equation} \label{Jutila-conj2}
\int_2^{T}\left(E(t+U)-E(t)\right)^4dx \ll T^{1+\varepsilon}U^2
\end{equation}
are true for $1\ll U\ll T^{1/2}$. He also pointed out that if \eqref{Jutila-conj2}
were true, then the important bound
\begin{equation*}
\int_1^T\left|\zeta(\frac{1}{2}+it)\right|^6dt\ll T^{1+\varepsilon}
\end{equation*}
would follow.

 When $H=T$, Recently Ivi\'c\cite{I2} obtained substantial improvements on  this
problem. He proved a much more explicit asymptotic formula for the
integral $\int_T^{2T}\left(\Delta(x+U)-\Delta(x)\right)^2dx.$ He
also proved that the estimate  (1.11)   holds for $T^{3/8}\ll U \ll
T^{1/2}.$
 Similar results
were also established for $E(t).$ We remark that the range
$T^{3/8}\ll U \ll T^{1/2}$ seems to be the limit of the present
methods.

Kiuchi and Tanigawa studied the mean square for short intervals of
Jutila's type for other arithmetical functions \cite{KT1, KT2}.

\medskip

\noindent {\bf Notation.} $\varepsilon$ is a sufficiently small positive
constant. $f\ll g$ means $|f|\leq Cg$ for some positive constant
$C.$ $n \asymp N$ means $N\ll n \ll N,$ $n\sim N$ means $N<n\leq
2N.$ $\mu(n)$ is the M\"obius function. $SC(\Sigma)$ denotes the
summation conditions of the sum $\Sigma$ when it is complicated.

\bigskip

%********************  section 2  *****************************************

\section{Main results}

 Jutila and Ivi\'c established asymptotic formulas for the mean
 square of
$ \Delta(x+U)-\Delta(x)  .$ However for the fourth moment, only an
upper bound result was proved, which seems very weak comparison to
the mean square case.  So it is natural to ask if one can find an
asymptotic formula for the fourth power moment of $
\Delta(x+U)-\Delta(x) $   in some range of $U.$ The aim of this
paper is to solve this problem . We shall establish an asymptotic
formula for the fourth power moment of $\Delta(x+U)-\Delta(x),$
which can be viewed as an analogue of (1.8).

\begin{thm}
Suppose $T,H,U$ are large real numbers such that
\begin{equation} \label{jyouken-1}
T^{3/7} \ll U\ll T^{1/2},\ \ H\leq T,  \ \ H^8U^{21}\gg T^{17},
\end{equation}
then for a small constant $c$ we have
\begin{align}
& \int_T^{T+H}(\Delta(x+U)-\Delta(x))^4dx \\
& =\frac{3}{2\pi^4} \!\!\!\! \sum_{\begin{subarray}{c} n_j \le c(T/U)^{1/4} \\
                                   \sqrt{n_1}+\sqrt{n_2}=\sqrt{n_3}+\sqrt{n_4}
       \end{subarray}} \!\!\!\! \frac{d(n_1)\cdots d(n_4)}{(n_1\cdots n_4)^{3/4}}
   \int_T^{T+H} x \prod_{j=1}^{4}\sin\frac{\pi U \sqrt{n_j}}{x^{1/2}}dx \nonumber \\
& \quad
+O\left(T^{17/8+\varepsilon}U^{-5/8}+HT^{3/16+\varepsilon}U^{25/16}
        +  H^{1/4}T^{15/16+\varepsilon}U^{21/16}\right) \nonumber
\end{align}
and
\begin{align}
& \int_T^{T+H}(\Delta(x+U)-\Delta(x))^4dx
=\frac{3}{\pi^4} \int_T^{T+H} x \left(\sum_{n\leq c(T/U)^{1/4}}\frac{d^2(n)}{n^{3/2}}
\sin\frac{\pi U \sqrt{n}}{x^{1/2}} \right)^2 dx \\
&\hspace{35mm}+O\left(HU^2\left(T/U^2\right)^{-1/2}T^{\varepsilon}+
HT^{3/16+\varepsilon}U^{25/16}\right) \nonumber\\
&\hspace{35mm}+O\left(T^{17/8+\varepsilon}U^{-5/8}+H^{1/4}T^{15/16+\varepsilon}U^{21/16}\right).\nonumber
\end{align}
\end{thm}

\noindent {\bf Remark 1.}  Our result holds  only in the range
$T^{3/7} \ll U \ll T^{1/2},$ which is   narrower than Ivi\'c's range
$T^{3/8} \ll U \ll T^{1/2}.$ However, our theorem gives an
asymptotic formula for the fourth power moment of
$\Delta(x+U)-\Delta(x).$ We note that the range $T^{3/7} \ll U \ll
T^{1/2}$ seems to be the limit  we can establish asymptotic result
for the fourth power moment by the present methods.

\noindent {\bf Remark 2.} If the condition (2.1) is replaced by
$$T^{3/7} \ll U\ll T^{1/2},\ \ T^{5/6}\ll H\leq T,  \ \ H^8U^{21}\gg T^{17},$$
then the term $H^{1/4}T^{15/16+\varepsilon}U^{21/16}$ in (2.2) and
(2.3) can be removed. If (2.1) is replaced by
$$T^{3/7} \ll
U\ll T^{1/2},\ \  H\leq T,  \ \ H^{16}U^{36}\gg T^{31},$$ then both
$H^{1/4}T^{15/16+\varepsilon}U^{21/16}$ and
$T^{17/8+\varepsilon}U^{-5/8}$ in (2.2) and (2.3) can be removed.

\begin{cor}
Suppose
\begin{equation} \label{jyouken-2}
 T^{3/7}\ll U \ll T^{1/2-\varepsilon}, \ \ H\leq T,\ \ H^{8}U^{21} \gg T^{17+\varepsilon}.
\end{equation}
Then we have
\begin{equation}
\int_T^{T+H}(\Delta(x+U)-\Delta(x))^4dx\asymp HU^2\log^6
\frac{T}{U^2}.
\end{equation}
\end{cor}

\smallskip

\smallskip

\begin{thm} Suppose $T,H,U$ are large real numbers such that
$$T^{3/7} \ll U\ll T^{1/2},\ \ T^{205/227}\ll H\leq T,  \ \ H^8U^{21}\gg T^{17},
$$
$c$ is a  small positive constant, then we have
\begin{align}
&\int_{T}^{T+H}\left(E(t+U)-E(t)\right)^4dt=\frac{12}{\pi}\sum_{\stackrel{\sqrt{n_1}+\sqrt{n_2}
=\sqrt{n_3}+\sqrt{n_4}}{n_j\leq
c(T/U)^{1/4}}}\frac{d(n_1)d(n_2)d(n_3)d(n_4)}{(n_1n_2n_3n_4)^{3/4}}\\
& \hspace{10mm}\times \int_T^{T+H}t\prod_{j=1}^4 \sin
\frac{U\sqrt{\pi n_j}}{(2t)^{1/2}}dt
+O\left(HT^{3/16+\varepsilon}U^{25/16}+T^{17/8+\varepsilon}U^{-5/8}\right)\nonumber
\end{align}
and
\begin{align}
&\int_T^{T+H}(E(t+U)-E(x))^4dt\\
&=\frac{24}{\pi}\sum_{n\leq c(T/U)^{1/4}}\frac{d^2(n)}{n^{3/2}}\sum_{m\leq c(T/U)^{1/4}}
\frac{d^2(m)}{m^{3/2}} \int_T^{T+H}t\sin^2\frac{U\sqrt {n\pi}}{(2t)^{1/2}}\sin^2
\frac{U\sqrt{m\pi}}{(2t)^{1/2}}dt\nonumber \\
&\hspace{20mm}+O\left(HU^2(\frac{T}{U^2})^{-1/2}T^{\varepsilon}
+HT^{3/16+\varepsilon}U^{25/16}+T^{17/8+\varepsilon}U^{-5/8}
\right).\nonumber
\end{align}
\end{thm}

\begin{cor} Suppose
$$T^{3/7} \ll U\ll T^{1/2-\varepsilon},\ \ T^{205/227}\ll H\leq T,  \ \ H^8U^{21}\gg T^{17+\varepsilon},
$$ then
\begin{equation}
\int_T^{T+H}(E(x+U)-E(x))^4dx\asymp HU^2\log^6 \frac{T}{U^2}.
\end{equation}
\end{cor}

\smallskip

%********************  section 3  **************************************

\section{Some preliminary lemmas}

\begin{lem} If $1\ll N\ll x,$ then
$$
\Delta(x)=\frac{x^{1/4}}{\pi\sqrt 2}\sum_{n\leq
N}\frac{d(n)}{n^{3/4}}\cos\left(4\pi\sqrt{nx}-\frac{\pi}{4}\right)
+O(x^{1/2+\varepsilon}N^{-1/2}).
$$
\end{lem}

\begin{proof}
This is the well-known truncated Voronoi's formula (see for example, (3.17) of Ivi\'c \cite{I}).
\end{proof}

\smallskip

\begin{lem} If $\alpha^{*}=\sqrt{n_1}+\sqrt{n_2}\pm\sqrt{n_3}-
\sqrt{n_4}\not=0$, then
$$
|\alpha^{*}| \gg \max(n_1,n_2,n_3, n_4)^{-3/2}(n_1n_2n_3n_4)^{-1/2}.
$$
\end{lem}

\begin{proof}
This is a variant of Lemma 2 of Ivi\'c and Sargos \cite{IS}.
\end{proof}

\smallskip

\begin{lem}
Let $N\geq 2$, $\Delta>0$ and let ${\cal A}(N;\Delta)$ denote the
number of solutions of the inequality
$$
|n_1^{1/2}+n_2^{1/2}- n_3^{1/2}-n_4^{1/2}|<\Delta,\quad  n_j\sim N \
(j=1,2,3,4),
$$
then
$$
{\cal A}(N;\Delta) \ll  (\Delta N^{7/2}+N^{2}) N^\varepsilon.
$$
\end{lem}

\begin{proof}
This is a special case of Theorem 2 of  Robert and Sargos \cite{RS}.
\end{proof}

\smallskip

\begin{lem}
Let $N_j \geq 2$, $\Delta>0$ and let ${\cal
A}_{\pm}(N_1,N_2,N_3,N_4;\Delta)$ denote the number of solutions of
the inequality
$$
|n_1^{1/2}+n_2^{1/2} \pm n_3^{1/2}-n_4^{1/2}|<\Delta,\quad  n_j\sim N \ (j=1,2,3,4),
$$
then
$$
{\cal A}_{\pm}(N_1,N_2,N_3,N_4;\Delta) \ll \prod_{j=1}^4 (\Delta^{1/4} N_j^{7/8}+N_j^{1/2})
N_j^\varepsilon.
$$
\end{lem}

\begin{proof}
This is Lemma 3 of Zhai \cite{Z3}.
\end{proof}

\smallskip

\begin{lem} Suppose  $f_j(t)(1\leq j\leq k)$ and $g(t)$ are
 continuous , monotonic real-valued functions on $[a,b]$
 and let $g(t)$ have a continuous , monotonic derivative on $[a,b].$
If $| f_j(t)|\leq A_j (1\leq j\leq k),
 |g^{\prime}(t)|\gg \Delta$ for any $t\in [a,b],$ then
$$\int_a^bf_1(t)\cdots f_k(t)e(g(t))dt\ll A_1\cdots A_k\Delta^{-1}.$$

\end{lem}

\begin{proof}This is Lemma 15.3 of Ivi\'c \cite{I}. We note that  Lemma 3.5 is still true
if  the function $e(t)$ is replaced by $\cos t$ and $\sin t.$
\end{proof}

\smallskip

\begin{lem}
The square roots of different square-free numbers
are linearly independent over the integers.
\end{lem}

\begin{proof}
This is a classical result of Besicovitch \cite{Be}.
\end{proof}

%\begin{lem}
%Suppose $N\geq 2,$ $\Delta>0.$ Let ${\cal A}(N;\Delta)$ denote the
%number of solutions of the inequality $$|n_1^{1/2}+n_2^{1/2}-
%n_3^{1/2}-n_4^{1/2}|<\Delta,\hspace{5mm}n_j\sim N(j=1,2,3,4),$$ then
%$${\cal A}(N;\Delta)\ll (\Delta N^{7/2}+N^{2})N^\varepsilon.$$
%\end{lem}

%\begin{proof}
%This is a special case of Lemma 3 of Zhai \cite{Z3}.
%\end{proof}

%******************  section 4  **********************************

\section{Estimates of some special sums}

In this section we shall prove the  estimates of several special sums used in our proof.

\begin{lem} Let $z\geq 10.$ Then
\begin{align*}
&\sum_{n\leq z}d^2(n)n^{-1/2}\asymp z^{1/2}\log^3 z,\\
\intertext{and}
&\sum_{n>z}d^2(n)n^{-3/2}\asymp z^{-1/2}\log^3 z.
\end{align*}
\end{lem}

\begin{proof}
These two estimates follow from the partial summation formula and the well-known estimate
$$
\sum_{n\leq z}d^2(n)\asymp z\log^3 z.
$$
\end{proof}

\begin{lem} Let $z\geq 10.$ Define
\begin{eqnarray*}
c_{1}(z):=\sum_{\begin{subarray}{c} \sqrt{n_1}+\sqrt{n_2}+\sqrt{n_3}=\sqrt{n_4} \\
                                    n_j\leq z \end{subarray}}
          \frac{d(n_1)d(n_2)d(n_3)d(n_4)}{{\mathstrut (n_1n_2n_3)}^{3/4} {\mathstrut n_4}^{1/4}}.
\end{eqnarray*}
Then $c_{1}(z) \ll 1.$
\end{lem}

\begin{proof}
Suppose $\sqrt{n_1}+\sqrt{n_2}+ \sqrt{n_3}=\sqrt{n_4}.$ By Lemma
3.6, we have
\begin{equation*}
n_j=lm_j^2, \, m_1+m_2+m_3=m_4, \, \mu(l) \not= 0.
\end{equation*}
So we have (noting that $d(ab)\leq d(a)d(b)$)
\begin{align*}
c_{1}(z)
&\ll \sum_{\begin{subarray}{c} m_1+m_2+m_3=m_4 \\ lm_4^2\leq z \end{subarray}}
     \frac{d(lm_1^2)d(lm_2^2)d(lm_3^2)d(lm_4^2)}{{\mathstrut l}^{5/2} {\mathstrut (m_1m_2m_3)}^{3/2}
           {\mathstrut m_4}^{1/2}}\\
&\ll \sum_{l\leq z}\frac{d^4(l)}{l^{5/2}}
     \sum_{\begin{subarray}{c} m_1+m_2+m_3=m_4 \\  m_4\leq (z/l)^{1/2} \end{subarray}}
     \frac{d(m_1^2)d(m_2^2)d(m_3^2)d(m_4^2)}{{\mathstrut (m_1m_2m_3)}^{3/2} {\mathstrut m_4}^{1/2}}\\
&\ll \sum_{l\leq z}\frac{d^4(l)}{l^{5/2}}
     \sum_{m_1, m_2, m_3\leq (z/l)^{1/2}}
     \frac{d(m_1^2)d(m_2^2)d(m_3^2)d((m_1+m_2+m_3)^2)}{(m_1m_2m_3)^{3/2}(m_1+m_2+m_3)^{1/2}}\\
& \ll 1.
\end{align*}
\end{proof}

\begin{lem} Let $z\geq 10.$ Define
\begin{eqnarray*}
c_{2}(z):=\sum_{\begin{subarray}{c} \sqrt{n_1}+\sqrt{n_2}=\sqrt{n_3}+\sqrt{n_4}\\
                                    n_j\leq z \end{subarray}}
          \frac{d(n_1)d(n_2)d(n_3)d(n_4)(\sqrt{n_1}+\sqrt{n_2}+\sqrt{n_3}+\sqrt{n_4})}{(n_1n_2n_3n_4)^{3/4}}.
\end{eqnarray*}
Then $c_{2}(z)\ll \log^4 z.$
\end{lem}

\begin{proof}
We have $c_{2}(z)\ll c_{2}^{\prime}(z),$ where
\begin{equation*}
c_{2}^{\prime}(z):=\sum_{\begin{subarray}{c} \sqrt{n_1}+\sqrt{n_2}=\sqrt{n_3}+\sqrt{n_4} \\
                                             n_j\leq z \end{subarray}}
      \frac{d(n_1)d(n_2)d(n_3)d(n_4)}{{\mathstrut (n_1n_2n_3)}^{3/4}{\mathstrut n_4}^{1/4}}.
\end{equation*}

It is easily seen that the equation $\sqrt{n_1}+\sqrt{n_2} =\sqrt{n_3}+\sqrt{n_4}$ is equivalent to
the following two cases:

(i) $n_1=n_3, n_2=n_4$ or $n_1=n_4, n_2=n_3,$

(ii) $n_1\not= n_3, n_1\not= n_4.$  \\
Thus we may write
\begin{align*}
c_{2}^{\prime}(z)&=2c_{21}(z)+c_{22}(z)\\
\intertext{with}
c_{21}(z)&:=\sum_{n\leq z}\frac{d^2(n)}{n^{3/2}}\sum_{m\leq z}\frac{d^2(m)}{m},\\
c_{22}(z)&:=\mathop{\sum\nolimits^{*}}_{\begin{subarray}{c}
            \sqrt{n_1}+\sqrt{n_2}=\sqrt{n_3}+\sqrt{n_4}\\ n_j\leq z \end{subarray}}
            \frac{d(n_1)d(n_2)d(n_3)d(n_4)}{{\mathstrut (n_1n_2n_3)}^{3/4}{\mathstrut n_4}^{1/4}},
\end{align*}
where $\sum^{*}$ means the condition $n_1\not= n_3, n_1\not= n_4.$

Obviously, $c_{21}(z)\ll \log^4 z.$  Now suppose
$\sqrt{n_1}+\sqrt{n_2} =\sqrt{n_3}+\sqrt{n_4}$ such that $n_1\not=
n_3, n_1\not= n_4.$ By Lemma 3.6 we have
\begin{equation}
n_j=lm_j^2,\hspace{2mm} m_1+m_2=m_3+m_4, \hspace{2mm}\mu(l)\not=0.
\end{equation}
Thus
\begin{align*}
c_{22}(z)
&\ll \sum_{l\leq z}\frac{d^4(l)}{l^{10/4}}
\sum_{\begin{subarray}{c} m_1+m_2=m_3+m_4 \\ m_j\leq (z/l)^{1/2} \end{subarray}}
\frac{d(m_1^2)d(m_2^2)d(m_3^2)d(m_4^2)}
{{\mathstrut (m_1m_2m_3)}^{3/2} {\mathstrut m_4}^{1/2}}\\
&\ll \sum_{l\leq z}\frac{d^4(l)}{l^{10/4}}\sum_{n\leq (z/l)^{1/2}}r(n)R(n),
\end{align*}
where
$$
r(n):=\sum_{n=m_1+m_2}\frac{d(m_1^2)d(m_2^2)}{(m_1m_2)^{3/2}} \quad \text{and} \quad
R(n)=\sum_{n=m_3+m_4}\frac{d(m_3^2)d(m_4^2)}{m_3^{3/2}m_4^{1/2}}.
$$
Obviously
$$
r(n)\ll n^{-3/2+\varepsilon} \quad \text{and} \quad R(n)\ll n^{-1/2+\varepsilon}.
$$
So we have $c_{22}(z)\ll 1.$ Now Lemma 4.3 follows from the above estimates.
\end{proof}

\begin{lem} Suppose $u>z\geq 10.$ Define
\begin{align*}
& c(z,u)=\sum\nolimits_1\frac{1}{(n_1n_2n_3n_4)^{3/4}} \prod_{j=1}^4
  \min\left(\frac{\sqrt n_j}{\sqrt z},1\right), \\
& SC(\Sigma_1): \sqrt{n_1}+\sqrt{n_2}=\sqrt{n_3}+\sqrt{n_4}, \, n_1\not= n_3, \, n_1\not= n_4,\,
   n_j\leq u \ \ (j=1,2,3,4).
\end{align*}
Then we have $c(z,u)\ll z^{-3/2}.$
\end{lem}

\begin{proof} Let
$$
g(n_1,n_2,n_3,n_4)=\frac{1}{(n_1n_2n_3n_4)^{3/4}}
\prod_{j=1}^4\min\left(\frac{\sqrt n_j}{\sqrt z},1\right).
$$
Obviously we have
\begin{align*}
& c(z,u)\ll \sum\nolimits_2g(n_1,n_2,n_3,n_4), \\
& SC(\Sigma_2): \sqrt{n_1}+\sqrt{n_2}=\sqrt{n_3}+\sqrt{n_4},\, n_1\not= n_3, \, n_1\not= n_4, \,
   n_j\leq u \ \ (j=1,2,3,4),\\
& \hspace{18mm}  n_1\leq n_2, \, n_3\leq n_4, \, n_2<n_4.
\end{align*}

We write
\begin{align*}
&\sum\nolimits_2 g(n_1,n_2,n_3,n_4)=(\sum\nolimits_3+\sum\nolimits_4)g(n_1,n_2,n_3,n_4), \\[1ex]
&SC(\Sigma_3): \sqrt{n_1}+\sqrt{n_2}=\sqrt{n_3}+\sqrt{n_4},\, n_1\not= n_3,\, n_1\not= n_4,\\
& \hspace{18mm} n_1\leq n_2, \, n_3\leq n_4, \, n_2<n_4, \, n_4\leq z,\\
&SC(\Sigma_4): \sqrt{n_1}+\sqrt{n_2}=\sqrt{n_3}+\sqrt{n_4},n_1\not= n_3, \, n_1\not= n_4,
               \, n_j\leq u \ \ (j=1,2,3,4),\\
& \hspace{18mm} n_1\leq n_2, \, n_3\leq n_4, \, n_2\leq n_4, \, n_4>z.
\end{align*}

We estimate $\Sigma_3$ first. From (4.1) we get
\begin{align*}
& \sum\nolimits_3 g(n_1,n_2,n_3,n_4) \ll z^{-2}\sum\nolimits_2 (n_1n_2n_3n_4)^{-1/4}\\
& \ll z^{-2}\sum_{l<z}l^{-1} \sum_{\begin{subarray}{c} m_1+m_2=m_3+m_4 \\
                                                       m_j\leq (z/l)^{1/2}
                                   \end{subarray}} (m_1m_2m_3m_4)^{-1/2}\\
& \ll z^{-2}\sum_{l<z}l^{-1}\sum_{n\leq 2(z/l)^{1/2}}f^2(n) \ll
      z^{-2}\sum_{l<z}l^{-1}\sum_{n\leq 2(z/l)^{1/2}}1\\
& \ll z^{-2}\sum_{l<z}l^{-1}(z/l)^{1/2}\ll z^{-3/2},
\end{align*}
where we used the estimate
$$
f(n)=\sum_{n=m_1+m_2}(m_1m_2)^{-1/2}\ll 1.
$$

Now we estimate $\Sigma_4.$ From the condition $SC(\Sigma_4)$ we
get $n_2\geq z/2.$ If $n_3\geq z,$ then $n_1\geq z.$ So we can
write
\begin{align*}
&\sum\nolimits_4 g(n_1,n_2,n_3,n_4)
 =\left(\sum\nolimits_{41}+\sum\nolimits_{42}+\sum\nolimits_{43}\right)g(n_1,n_2,n_3,n_4),\\[1ex]
&SC(\Sigma_{41}):\sqrt{n_1}+\sqrt{n_2}=\sqrt{n_3}+\sqrt{n_4},n_1\not= n_3, n_1\not= n_4,\\
&\hspace{18mm} n_4>z, n_2>z/2, n_3<z, n_1<z,\\
&SC(\Sigma_{42}):\sqrt{n_1}+\sqrt{n_2}=\sqrt{n_3}+\sqrt{n_4},n_1\not= n_3, n_1\not= n_4,\\
&\hspace{18mm} n_4>z, n_2>z, n_3<z, n_1>z,\\
&SC(\Sigma_{43}):\sqrt{n_1}+\sqrt{n_2}=\sqrt{n_3}+\sqrt{n_4},n_1\not= n_3, n_1\not= n_4,\\
&\hspace{18mm} n_4>z, n_2>z, n_3\geq z, n_1\geq z.\\
\end{align*}

By (4.1) we get
\begin{align*}
& \sum\nolimits_{41}g(n_1,n_2,n_3,n_4)
  \ll z^{-1}\sum\nolimits_{41}n_1^{-1/4}n_2^{-3/4}n_3^{-1/4}n_4^{-3/4}\\
& \ll z^{-1}\sum_{l}l^{-2}\sum\nolimits_{411}m_1^{-1/2}m_2^{-3/2}m_3^{-1/2}m_4^{-3/2}\\
& =z^{-1}\sum_{l}l^{-2}\sum_{n>(z/l)^{1/2}}f_1(n)g_1(n),
\end{align*}
where
\begin{align*}
& SC(\Sigma_{411}): m_1+m_2=m_3+m_4, m_1\leq (z/l)^{1/2}, m_3\leq (z/l)^{1/2},\\
& \hspace{18mm} m_2>(z/2l)^{1/2}, m_4>(z/l)^{1/2};\\
& f_1(n):=\sum\nolimits_{412}m_1^{-1/2}m_2^{-3/2},\hspace{2mm}
  SC(\Sigma_{412}):m_1\leq (z/l)^{1/2}, m_2>(z/2l)^{1/2},m_1\leq m_2;\\
& g_1(n):=\sum\nolimits_{413}m_3^{-1/2}m_4^{-3/2},\hspace{2mm}
  SC(\Sigma_{413}): m_3\leq (z/l)^{1/2}, m_4>(z/l)^{1/2}.
\end{align*}
It is easy to see that
$$
f_1(n)\ll n^{-3/2}(z/l)^{1/4},\hspace{2mm} g_1(n)\ll n^{-3/2}(z/l)^{1/4}.
$$
Thus we have
\begin{align*}
& \sum\nolimits_{41}g(n_1,n_2,n_3,n_4)\ll z^{-1}\sum_{l}l^{-2}(z/l)^{1/2}
\sum_{n>(z/l)^{1/2}}n^{-3}\ll z^{-3/2}.
\end{align*}
Similarly we can prove
\begin{align*}
& \sum\nolimits_{42}g(n_1,n_2,n_3,n_4)\ll z^{-3/2},\\
& \sum\nolimits_{43}g(n_1,n_2,n_3,n_4)\ll z^{-3/2}.
\end{align*}
Now Lemma 4.4 follows from the above estimates.
\end{proof}

%*************************** section 5 **********************************************

\section{\bf Identities involving the functions $\sin u$ and $\cos u$}

We need some formulas about the functions $\sin u$ and $\cos u.$ Let
$k\geq 1$ be a fixed integer, $\alpha_1,\cdots, \alpha_k$ be real
numbers. Define
\begin{align*}
\mathcal{SI}_k(\alpha_1,\cdots, \alpha_k)&:=\sum_{j_1=0}^1\cdots
\sum_{j_k=0}^1(-1)^{j_1+\cdots+j_k}\sin(j_1\alpha_1+\cdots+j_k\alpha_k),\\
\mathcal{CO}_k(\alpha_1,\cdots, \alpha_k)&:=\sum_{j_1=0}^1\cdots
\sum_{j_k=0}^1(-1)^{j_1+\cdots+j_k}\cos(j_1\alpha_1+\cdots+j_k\alpha_k).
\end{align*}

\begin{lem} Suppose $k=m+l,$  $m\geq 1, l\geq 1,$ then we have
\begin{align}
\mathcal{SI}_k(\alpha_1,\cdots, \alpha_k)
&=\mathcal{SI}_m(\alpha_1,\cdots, \alpha_m)\mathcal{CO}_l(\alpha_{m+1},\cdots, \alpha_k)\\
&\quad +\mathcal{CO}_m(\alpha_1,\cdots,\alpha_m)
        \mathcal{SI}_l(\alpha_{m+1},\cdots, \alpha_k),\nonumber\\
\mathcal{CO}_k(\alpha_1,\cdots, \alpha_k)
&=\mathcal{CO}_m(\alpha_1,\cdots, \alpha_m)\mathcal{CO}_l(\alpha_{m+1},\cdots, \alpha_k)\\
&\quad -\mathcal{SI}_m(\alpha_1,\cdots,\alpha_m)\mathcal{SI}_l(\alpha_{m+1},\cdots, \alpha_k).\nonumber
\end{align}
\end{lem}

\begin{proof}
These two identities follow from the well-known  formulas
$$\sin(\alpha+\beta)=\sin\alpha\cos\beta+\cos\alpha\sin\beta$$ and
 $$\cos(\alpha+\beta)=\cos\alpha\cos\beta-\sin\alpha\sin\beta,$$
respectively.
\end{proof}

\begin{lem} We have
\begin{align*}
\mathcal{SI}_k(\alpha_1,\cdots, \alpha_k)=
\left\{\begin{array}{l}
2^k(-1)^{\frac{r+1}{2}}\left(\prod_{j=1}^k\sin\frac{\alpha_j}{2}\right)
\cos\frac{\alpha_1+\cdots+\alpha_k}{2}, \\
\mbox{\hspace{4cm} if \ $k\equiv r \pmod{4}, \  r=1,3,$}\\
2^k(-1)^{\frac{r}{2}}\left(\prod_{j=1}^k\sin\frac{\alpha_j}{2}\right)
\sin\frac{\alpha_1+\cdots+\alpha_k}{2}, \\
\mbox{\hspace{4cm} if \ $k\equiv r \pmod{4}, \ r=2,4,$}
\end{array}\right.
\end{align*}
\begin{align*}
\mathcal{CO}_k(\alpha_1,\cdots, \alpha_k)=\left\{\begin{array}{l}
2^k(-1)^{\frac{r-1}{2}}\left(\prod_{j=1}^k\sin\frac{\alpha_j}{2}\right)
\sin\frac{\alpha_1+\cdots+\alpha_k}{2}, \\
\mbox{\hspace{4cm} if \ $k\equiv r \pmod{4}, r=1,3,$}\\
2^k(-1)^{\frac{r}{2}}\left(\prod_{j=1}^k\sin\frac{\alpha_j}{2}\right)
\cos\frac{\alpha_1+\cdots+\alpha_k}{2}, \\
\mbox{\hspace{4cm} if \ $k\equiv r \pmod{4}, r=2,4.$}
\end{array}\right.
\end{align*}
\end{lem}

{\allowdisplaybreaks

\begin{proof}
Trivially for $k=1$ we have
\begin{align}
\mathcal{SI}_1(\alpha_1)&=-\sin\alpha_1=-2\sin\frac{\alpha_1}{2}\cos\frac{\alpha_1}{2},\\
\mathcal{CO}_1(\alpha_1)&=1-\cos\alpha_1=2\sin^2\frac{\alpha_1}{2}. %\nonumber
\end{align}
The cases $k=2, 3, 4$ follow easily  from Lemma 5.1 and  (5.3).

Now suppose $m\geq 1$ be a fixed integer such that Lemma 5.2 is
true for $k\leq 4m.$ We shall show that Lemma 5.2 is true for
$k=4m+r,$ $r=1,2,3,4.$ We only prove the result for
$\mathcal{SI}_k(\alpha_1,\cdots, \alpha_k)$ with $r=1,3.$ The other cases
are similar. By Lemma 5.1 we have
\begin{align*}
&\mathcal{SI}_k(\alpha_1,\cdots, \alpha_k)=\mathcal{SI}_{4m}(\alpha_1,\cdots,
\alpha_{4m}) \mathcal{CO}_r(\alpha_{4m+1},\cdots,\alpha_k)\\
&\hspace{3cm}  + \mathcal{CO}_{4m}(\alpha_1,\cdots,\alpha_{4m})
\mathcal{SI}_r(\alpha_{4m+1},\cdots, \alpha_k)\\
&=(-1)^{\frac{r-1}{2}}2^k\left\{\prod_{j=1}^k\sin\frac{\alpha_j}{2}\right\}
\sin\frac{\alpha_1+\cdots+\alpha_{4m}}{2}\sin\frac{\alpha_{4m+1}+\cdots+\alpha_k}{2}\\
&\hspace{8mm}+(-1)^{\frac{r+1}{2}}2^k\left\{\prod_{j=1}^k\sin\frac{\alpha_j}{2}\right\}
\cos\frac{\alpha_1+\cdots+\alpha_{4m}}{2}\cos\frac{\alpha_{4m+1}+\cdots+\alpha_k}{2}\\
&=(-1)^{\frac{r+1}{2}}2^k\left\{\prod_{j=1}^k\sin\frac{\alpha_j}{2}\right\}
\{\cos\frac{\alpha_1+\cdots+\alpha_{4m}}{2}
\cos\frac{\alpha_{4m+1}+\cdots+\alpha_k}{2}\\
&\hspace{30mm}-\sin\frac{\alpha_1+\cdots+\alpha_{4m}}{2}
\sin\frac{\alpha_{4m+1}+\cdots+\alpha_k}{2}\}\\
&=(-1)^{\frac{r+1}{2}}2^k\left\{\prod_{j=1}^k\sin\frac{\alpha_j}{2}\right\}
\cos\frac{\alpha_{1}+\cdots+\alpha_k}{2}.
\end{align*}

\end{proof}
}

%************************* section 6 *********************************************

\section{\bf Proof of Theorem 1}

Suppose the condition \eqref{jyouken-1} holds. Let  $y:= c(T/U)^{1/4}$, where
$c$ is a small positive constant. For any $T\leq x\ll T$, define
\begin{align*}
&R_1(x):=\frac{x^{1/4}}{\pi \sqrt{2}}\sum_{n\leq y}\frac{d(n)}{n^{3/4}}
         \cos(4\pi\sqrt{xn}-\frac{\pi}{4}),\\
&R_2(x):=\Delta(x)-R_1(x).
\end{align*}
By Lemma 3.1, we can write
\begin{equation}
\Delta(x+U)-\Delta(x)=S_1(x)+S_2(x),
\end{equation}
where
\begin{align*}
S_1(x)&=\frac{x^{1/4}}{\pi\sqrt 2}\sum_{n\leq y}\frac{d(n)}{n^{3/4}}
\left\{\cos\left(4\pi\sqrt{n(x+U)}-\frac{\pi}{4}\right)-
\cos\left(4\pi\sqrt{nx}-\frac{\pi}{4}\right)\right\},\\
S_2(x)&=R_2(x+U)-R_2(x)+M(x),
\end{align*}
and
\begin{align*}
M(x)&=\frac{(x+U)^{1/4}-x^{1/4}}{\pi\sqrt 2}\sum_{n\leq y}\frac{d(n)}{n^{3/4}}
\cos\left(4\pi\sqrt{n(x+U)}-\frac{\pi}{4}\right) \\
&=\frac{(x+U)^{1/4}-x^{1/4}}{(x+U)^{1/4}}{\cal R}_1(x+U).
\end{align*}
It is easy to see that
$$
M(x) \ll Ux^{-1}|R_1(x+U)|.
$$

%**************** subsection 6.1  **************************************

\subsection{Evaluation of the integral $\int_T^{T+H}S_1^4(x)dx$}

First we evaluate the integral $\int_T^{T+H}S_1^4(x)dx$. $S_1(x)$ can be written as
$$
S_1(x)=\frac{x^{1/4}}{\pi\sqrt 2}\sum_{j=0}^1(-1)^{j+1}\sum_{n\leq y}\frac{d(n)}{n^{3/4}}
\cos\left(4\pi\sqrt{n(x+jU)}-\frac{\pi}{4}\right).
$$
In order to write the 4-th power in a simple way, we introduce here the following
notations. Let $I=\{0,1\}$ and let $\bvec{j}=(j_1,j_2,j_3,j_4)$ and $\bvec{i}=(i_1,i_2,i_3)$
be elements in $I^4$ and $I^3$, respectively, and put $|\bvec{j}|=j_1+\cdots+j_4$.

Using the elementary formula
$$
\cos a_1\cos a_2\cos a_3 \cos a_4= \frac{1}{2^{3}}\sum_{\bvec{\scriptstyle i}\in I^{3}}
\cos(a_1+(-1)^{i_1}a_2+(-1)^{i_2}a_3+(-1)^{i_{3}}a_4)
$$
we have
\begin{align} \label{eq62}
S_1^4(x)&=\frac{x}{(\pi\sqrt 2)^4}\sum_{\bvec{\scriptstyle j}\in I^4}(-1)^{|\bvec{\scriptstyle j}|+4}
\sum_{n_1\leq y}\cdots\sum_{n_4\leq y}\frac{d(n_1)\cdots d(n_4)}{(n_1\cdots n_4)^{3/4}}\\
&\hspace{30mm}\times \prod_{l=1}^4 \cos\left(4\pi\sqrt{n_l(x+j_lU)}-\frac{\pi}{4}\right)\nonumber\\
&=\frac{x}{2^5\pi^4}\sum_{\bvec{\scriptstyle j}\in
I^4}(-1)^{|\bvec{\scriptstyle j}|}\sum_{\bvec{\scriptstyle i}\in I^3}
\sum_{n_i \leq y} \frac{d(n_1)\cdots d(n_4)}{(n_1\cdots n_4)^{3/4}}
%\sum_{n_1\leq y}\cdots\sum_{n_4\leq y}\frac{d(n_1)\cdots d(n_4)}{(n_1\cdots n_4)^{3/4}}
\cos\left(4\pi\alpha(x)-\frac{\pi\beta}{4}\right),\nonumber
\end{align}
where
\begin{align*}
\alpha(x)&=\sqrt{n_1(x+j_1U)}+(-1)^{i_1}\sqrt{n_2(x+j_2U)}+(-1)^{i_2}\sqrt{n_3(x+j_3U)}\\
         & \quad + (-1)^{i_{3}}\sqrt{n_4(x+j_4U)},\\
\beta &= 1+(-1)^{i_1}+(-1)^{i_2}+(-1)^{i_{3}}.
\end{align*}
Let
\begin{align}
&\alpha^{*}=\sqrt{n_1}+(-1)^{i_1}\sqrt{n_2}+(-1)^{i_2}\sqrt{n_3}+(-1)^{i_{3}}\sqrt {n_4},
\label{eq63} \\[1ex]
&\alpha_{*}=j_1\sqrt{n_1}+(-1)^{i_1}j_2\sqrt{n_2}+(-1)^{i_2}j_3\sqrt{n_3}+(-1)^{i_{3}}j_4
\sqrt {n_4}.  \label{eq64}
\end{align}
We divide the sum \eqref{eq62} into two parts according to $\alpha^*=0$ or
$\alpha^* \neq 0$, respectively. Thus we get
\begin{equation} \label{eq65}
S_1^4(x)=\frac{1}{ 2^5\pi^4}(S_{11}(x)+S_{12}(x)),
\end{equation}
where
\begin{align*}
S_{11}(x)
&=x\sum_{\bvec{\scriptstyle j} \in I^4} (-1)^{|\bvec{\scriptstyle j}|}
\sum_{\bvec{\scriptstyle i}\in I^3}
\sum_{\begin{subarray}{c} n_i \leq y \\ \alpha^{*}=0 \end{subarray}}
\frac{d(n_1)\cdots d(n_4)}{(n_1\cdots n_4)^{3/4}}\cos\left(4\pi\alpha(x)-\frac{\pi\beta}{4}\right),\\
\intertext{and}
S_{12}(x)
&=x\sum_{\bvec{\scriptstyle j} \in I^4} (-1)^{|\bvec{\scriptstyle j}|}
\sum_{\bvec{\scriptstyle i}\in I^3}
  \sum_{\begin{subarray}{c} n_i \leq y \\ \alpha^{*} \neq 0 \end{subarray}}
  \frac{d(n_1)\cdots d(n_4)}{(n_1\cdots
  n_4)^{3/4}}\cos\left(4\pi\alpha(x)-\frac{\pi\beta}{4}\right).
\end{align*}

It is easy to see that
$$
4\pi\alpha(x)=4\pi\alpha^{*}\sqrt x+\frac{2\pi U\alpha_{*}}{x^{1/2}}
+O\left(\frac{U^2(\sqrt{n_1}+\cdots+\sqrt{n_4})}{x^{3/2}}\right).
$$
Hence, if $\alpha^{*}=0$, we have
\begin{align*}
\cos\left(4\pi\alpha(x)-\frac{\pi\beta}{4}\right)
&=\cos\left(\frac{2\pi U \alpha_{*}}{x^{1/2}}-\frac{\pi\beta}{4}
+O\left(\frac{U^2(\sqrt{n_1}+\cdots+\sqrt{n_4})}{x^{3/2}}\right)\right)\\
&=\cos\left(\frac{2\pi U \alpha_{*}}{x^{1/2}}-\frac{\pi\beta}{4}
\right)+O\left(\frac{U^2(\sqrt{n_1}+\cdots+\sqrt{n_4})}{x^{3/2}}\right).
\end{align*}
Therefore
\begin{align} \label{eq66}
&\int_T^{T+H}S_{11}(x)dx=\sum_{\bvec{\scriptstyle j} \in I^4}
(-1)^{|\bvec{\scriptstyle j}|} \sum_{\bvec{\scriptstyle i}\in I^3}
\sum_{\begin{subarray}{c} n_i \leq y \\ \alpha^{*} = 0
\end{subarray}} \frac{d(n_1)\cdots d(n_4)}{(n_1\cdots n_4)^{3/4}}\\
& \hspace{4cm}\times \int_T^{T+H}x\cos\left(\frac{2\pi U\alpha_{*}}{x^{1/2}}
-\frac{\pi\beta}{4} \right)dx    \nonumber\\
&\qquad \qquad{}+O\left(HT^{-1/2}U^2\sum_{\alpha^{*}=0}\frac{d(n_1)\cdots
d(n_4)(\sqrt{n_1}+\cdots+\sqrt{n_4})}{(n_1\cdots
n_4)^{3/4}}\right).\nonumber
\end{align}
For $\bvec{i}=(i_1,i_2,i_3)\in I^3, \bvec{i} \not= (0,0,0)$, we let
\begin{align*}
&I(\bvec{i})=\sum_{\bvec{\scriptstyle j} \in I^4} (-1)^{|\bvec{\scriptstyle j}|}
 \sum_{\begin{subarray}{c} n_i \leq y \\ \alpha^{*} = 0 \end{subarray}}
 \frac{d(n_1)d(n_2)d(n_3)d(n_4)}{(n_1n_2n_3n_4)^{3/4}}\\
& \hspace{5cm} \times \int_T^{T+H}x \cos\left(\frac{2\pi U \alpha_*}{x^{1/2}}
-\frac{\pi \beta}{4}\right)dx,\\
% &&SC(\Sigma_5):
%\sqrt{n_1}+(-1)^{i_1}\sqrt{n_2}+(-1)^{i_2}\sqrt{n_3}+(-1)^{i_3}\sqrt{n_4}=
%0, n_j\leq y(j=1,2,3,4).
\end{align*}
where $\alpha^*$ and $\alpha_*$ are defined by \eqref{eq63} and \eqref{eq64}.
It is easily seen that
\begin{gather*}
I(0,0,1)=I(0,1,0)=I(1,0,0)=I(1,1,1),\\
I(0,1,1)=I(1,1,0)=I(1,0,1).
\end{gather*}
Concerning the sum in $O$-term in \eqref{eq66}, Lemma 4.2 and 4.3
imply that
\begin{equation*}
\sum_{\begin{subarray}{c} n_i \leq y \\ \alpha^{*}=0 \end{subarray}}
\frac{d(n_1)\cdots d(n_4)(\sqrt{n_1}+\cdots+\sqrt{n_4})}{(n_1\cdots n_4)^{3/4}}\ll \log^4 y.
\end{equation*}
Therefore
\begin{align*}
&\int_{T}^{T+H}S_{11}(x)dx
=4\sum_{\begin{subarray}{c} \sqrt{n_1}+\sqrt{n_2}+\sqrt{n_3}=\sqrt{n_4} \\  n_j\leq y \end{subarray}}
 \frac{d(n_1)d(n_2)d(n_3)d(n_4)}{(n_1n_2n_3n_4)^{3/4}} \\
&\hspace{1cm} \times \int_T^{T+H} x \mathcal{SI}_4\left(\frac{2\pi
U\sqrt{n_1}}{x^{1/2}},
   \frac{2\pi U\sqrt{n_2}}{x^{1/2}},\frac{2\pi U\sqrt{n_3}}{x^{1/2}},
   -\frac{2\pi U\sqrt{n_4}}{x^{1/2}}\right)dx \\[1ex]
& \quad {}+3\sum_{\begin{subarray}{c} \sqrt{n_1}+\sqrt{n_2}=\sqrt{n_3}+\sqrt{n_4} \\ n_j\leq y
                  \end{subarray}}
            \frac{d(n_1)d(n_2)d(n_3)d(n_4)}{(n_1n_2n_3n_4)^{3/4}} \\
&\hspace{1cm} \times \int_T^{T+H}x \mathcal{CO}_4\left(\frac{2\pi
U\sqrt{n_1}}{x^{1/2}},
              \frac{2\pi U\sqrt{n_2}}{x^{1/2}},-\frac{2\pi U\sqrt{n_3}}{x^{1/2}},
             -\frac{2\pi U\sqrt{n_4}}{x^{1/2}}\right)dx\\[2ex]
&\quad  {}+O(HU^2T^{-1/2}\log^4 y).
\end{align*}
By Lemma 5.2, we have
$$
\mathcal{SI}_4\left(\frac{2\pi U\sqrt{n_1}}{x^{1/2}},\frac{2\pi U\sqrt{n_2}}{x^{1/2}},
\frac{2\pi U\sqrt{n_3}}{x^{1/2}},-\frac{2\pi U\sqrt{n_4}}{x^{1/2}}\right)=0$$
for $\sqrt{n_1}+\sqrt{n_2}+\sqrt{n_3}=\sqrt{n_4}$, and
$$
\mathcal{CO}_4\left(\frac{2\pi U\sqrt{n_1}}{x^{1/2}},\frac{2\pi U\sqrt{n_2}}{x^{1/2}},
-\frac{2\pi U\sqrt{n_3}}{x^{1/2}},-\frac{2\pi U\sqrt{n_4}}{x^{1/2}}\right)
=16\prod_{j=1}^4 \sin \frac{\pi U\sqrt{n_j}}{x^{1/2}}$$
for $\sqrt{n_1}+\sqrt{n_2}=\sqrt{n_3}+\sqrt{n_4}$.
Hence
\begin{align}  \label{eq67}
&\int_{T}^{T+H}S_{11}(x)dx
=48\sum_{\begin{subarray}{c} \sqrt{n_1}+\sqrt{n_2}=\sqrt{n_3}+\sqrt{n_4} \\
                            n_j\leq y \end{subarray}}
\frac{d(n_1)\cdots d(n_4)}{(n_1 \cdots n_4)^{3/4}}\\
&\hspace{5mm}\times \int_T^{T+H}x\prod_{j=1}^4 \sin \frac{\pi U\sqrt{n_j}}{x^{1/2}}dx
+O(HU^2T^{-1/2}\log^4 y).\nonumber
\end{align}

Now consider the integral $\int_{T}^{T+H}S_{12}(x)dx$.  Note that
$S_{12}(x)$ corresponds to the sum under the condition $\alpha^*
\neq 0$. From the definition of $\alpha(x)$, we have
$$
\alpha'(x)=\frac{\alpha^*}{2\sqrt{x}}
+O\left(\frac{U(\sqrt{n_1}+\cdots+\sqrt{n_4})}{x^{3/2}}\right).
$$
By Lemma 3.2 we get $|\alpha^{*}|\gg \max(n_1,n_2,n_3,n_4)^{-7/2},$
which combining with $y=c(T/U)^{1/4}$ with a small positive constant $c$ implies that
$$
\alpha^{\prime}(x)\gg |\alpha^{*}|T^{-1/2}.
$$
By Lemma 3.5  we get
$$
\int_{T}^{T+H}x\cos(4\pi \alpha(x)-\pi\beta/4)dx \ll \frac{T^{3/2}}{|\alpha^*|}.
$$
Hence
\begin{align*}
\int_T^{T+H}S_{12}(x)dx & \ll T^{3/2}\sum_{\bvec{\scriptstyle j}
\in I^{4}}\sum_{\bvec{\scriptstyle i} \in I^{3}}
\sum_{\begin{subarray}{c} n_i \leq y \\ \alpha^* \neq 0 \end{subarray}}
\frac{d(n_1)\cdots d(n_4)} {(n_1\cdots n_4)^{3/4}} \frac{1}{|\alpha^*|} \\
& \ll G_1+G_2,
\end{align*}
where
\begin{align*}
&G_1=T^{3/2} \sum_{\begin{subarray}{c} n_1,n_2,n_3,n_4\leq y \\
                     \sqrt{n_1}+\sqrt{n_2}\neq \sqrt{n_3}+\sqrt{n_4} \end{subarray}}
     \frac{d(n_1)\cdots d(n_4)} {(n_1\cdots n_4)^{3/4}}
     \frac{1}{|\sqrt{n_1}+\sqrt{n_2}-\sqrt{n_3}-\sqrt{n_4}|},\\
&G_2=T^{3/2} \sum_{\begin{subarray}{c} n_1,n_2,n_3,n_4\leq y \\
                     \sqrt{n_1}+\sqrt{n_2}+\sqrt{n_3} \neq \sqrt{n_4} \end{subarray}}
     \frac{d(n_1)\cdots d(n_4)} {(n_1\cdots n_4)^{3/4}}
     \frac{1}{|\sqrt{n_1}+\sqrt{n_2}+\sqrt{n_3}-\sqrt{n_4}|}.
\end{align*}

We only estimate the sum $G_1.$ The estimate for $G_2$ is the same.

Let $\eta=\sqrt{n_1}+\sqrt{n_2}-\sqrt{n_3}-\sqrt{n_4}$. By a splitting argument we get
\begin{equation*}
G_1 \ll G_1(N_1, N_2, N_3, N_4)\log^4 T
\end{equation*}
for some $N_1, N_2, N_3, N_4$,  $N_j \ll y$ ($j=1,2,3,4$), where
$$
G_1(N_1, N_2, N_3, N_4)=T^{3/2} \sum_{\begin{subarray}{c} n_i\sim N_i \\ \eta \neq 0 \end{subarray}}
\frac{d(n_1)\cdots d(n_4)} {(n_1\cdots n_4)^{3/4}} \frac{1}{| \eta|}.
$$
Without loss of generality, we may assume that $N1 \le N_2, N_3 \le N_4$ and $N_2 \le N_4$.
Then by Lemma 3.2 we have
\begin{equation} \label{eq68}
\eta_0:=N_4^{-3/2}(N_1N_2N_3N_4)^{-1/2} \ll |\eta| \ll N_4^{1/2}
\end{equation}
 By a splitting argument again we get that
\begin{equation*}
G_1(N_1, N_2, N_3, N_4)\ll
\frac{T^{3/2+\varepsilon}}{(N_1N_2N_3N_4)^{3/4}}\frac{1}{\delta}\sum_{\delta<|\eta|\leq 2\delta}1.
\end{equation*}
By Lemma 3.4 we have
\begin{align*}
&G_1(N_1, N_2, N_3, N_4) \\
& \ll \frac{T^{3/2+\varepsilon}}{(N_1N_2N_3N_4)^{3/4}\delta} \prod_{j=1}^4
      N_j^{1/2}(\delta^{1/4}N_j^{3/8}+1) \\
& \ll \frac{T^{3/2+\varepsilon}}{(N_1N_2N_3N_4)^{1/4}\delta}
      \left(\delta^{3/4}(N_1N_2N_3)^{3/8}+\delta^{1/2}N_4^{3/4}+\delta^{1/4}N_4^{3/8}+1\right)
       (\delta^{1/4}N_4^{3/8}+1) \\
& \ll \frac{T^{3/2+\varepsilon}}{(N_1N_2N_3N_4)^{1/4}\delta}
      \left(\delta^{3/4}(N_1N_2N_3)^{3/8}(\delta^{1/4}N_4^{3/8}+1)+\delta^{3/4}N_4^{9/8}+1\right)\\
& \ll \frac{T^{3/2+\varepsilon}(N_1N_2N_3)^{1/8}}{N_4^{1/4}}(N_4^{3/8}+\delta^{-1/4})
      +\frac{T^{3/2+\varepsilon}}{(N_1N_2N_3N_4)^{1/4}}(\delta^{-1/4}N_4^{9/8}+\delta^{-1}) \\
& \ll T^{3/2+\varepsilon}N_4^{5/2}
\end{align*}
where in the last step we used $(6.8).$ Thus we have $G_1\ll T^{3/2+\varepsilon}y^{5/2}$.
Similarly we have $G_2\ll T^{3/2+\varepsilon}y^{5/2}$.
Combining the above estimates we get
\begin{equation}  \label{eq69}
\int_T^{T+H}S_{12}(x)dx\ll T^{3/2+\varepsilon}y^{5/2}.
\end{equation}

From \eqref{eq65}, \eqref{eq67} and \eqref{eq69}, we get
\begin{align} \label{eq610}
\int_{T}^{T+H}S_{1}^4(x)dx&=\frac{3}{2\pi^4}
\sum_{\begin{subarray}{c} \sqrt{n_1}+\sqrt{n_2}=\sqrt{n_3}+\sqrt{n_4} \\ n_j\leq y \end{subarray}}
\frac{d(n_1)d(n_2)d(n_3)d(n_4)}{(n_1n_2n_3n_4)^{3/4}}\\
& \qquad  \times \int_T^{T+H}x\prod_{j=1}^4 \sin \frac{\pi U\sqrt{n_j}}{x^{1/2}}dx    \nonumber\\
& \quad +O(T^{3/2+\varepsilon}y^{5/2}+HU^2T^{-1/2}\log^4 T). \nonumber
\end{align}

It is easy to see that  the sum in the right hand side of \eqref{eq610} can
be  written as  $\sum_6+\sum_7$, where
\begin{align}
\sum\nolimits_6&=\frac{3}{\pi^4}\sum_{n\leq y}\frac{d^2(n)}{n^{3/2}}\sum_{m\leq y}\frac{d^2(m)}{m^{3/2}}
\int_T^{T+H}x\sin^2 \frac{\pi U\sqrt{n}}{x^{1/2}}\sin^2\frac{\pi U\sqrt{m}}{x^{1/2}}dx\\
&=\frac{3}{\pi^4}\int_{T}^{T+H}x\left(\sum_{n \le y} \frac{d^2(n)}{n^{3/2}}\sin^2
  \frac{\pi U \sqrt{n}}{x^{1/2}}\right)^2 dx, \nonumber\\
\intertext{and} \sum\nolimits_7&=\frac{3}{2\pi^4} \!\!\!
\sum_{\begin{subarray}{c} n_i \leq y \\ \sqrt{n_1}+\sqrt{n_2}
=\sqrt{n_3}+\sqrt{n_4}\\ n_1 \neq n_3, n_1 \neq n_4 \end{subarray}}
\!\!\! \frac{d(n_1)\cdots d(n_4)}{(n_1 \cdots n_4)^{3/4}}
\int_T^{T+H}x\prod_{j=1}^4 \sin \frac{\pi U\sqrt{n_j}}{x^{1/2}}dx.
\end{align}
By Lemma 4.1 we have
\begin{equation}
\sum\nolimits_6 \ll HU^2\log^6 T,
\end{equation}
while by Lemma 4.4 we have
\begin{equation} \label{sigma7}
\sum\nolimits_7 \ll T^{1+\varepsilon}H(T/U^2)^{-3/2} \ll HU^2
T^{\varepsilon}(T/U^2)^{-1/2}.
\end{equation}

From (6.10)--(6.14) we get
\begin{eqnarray}
\int_T^{T+H}S_1^4(x)dx\ll
HU^2T^{\varepsilon}+T^{3/2+\varepsilon}y^{5/2}\ll
HU^2T^{\varepsilon}
\end{eqnarray}
if we note the condition $H^8U^{21}\gg T^{17}$.

%************  subsection 6.2  ******************************************************

\subsection{On the integral $\int_T^{T+H}S_2^4(x)dx$}

In this subsection we estimate the integral $\int_T^{T+H}S_2^4(x)dx$.
Since $S_2(x)=R_2(x+U)-R_2(x)+M(x)$, it is sufficient to consider the integrals of
$R_2^4(x)$ and $M^4(x)$.

Recall that ($T\leq x\leq T+H$)
\begin{align*}
R_2(x)=\frac{x^{1/4}}{\pi \sqrt{2}}\sum_{y<n \le
T}\frac{d(n)}{n^{3/4}}\cos(4\pi \sqrt{nx}-\pi/4) +O(T^{\varepsilon}
).
\end{align*}

Let $J$ be a positive integer such that $2^{J+1}\asymp
H^2y^{-3}T^{-1}$  and let $z=2^{J+1}y$, then $y<z<T.$  We divide
$R_2(x)$ into two parts:
\begin{align*}
&R_2(x)=R_{21}(x)+R_{22}(x),\\
\intertext{where}
R_{21}(x)&=\frac{x^{1/4}}{\sqrt 2\pi}\sum_{y<n\leq z}
           \frac{d(n)}{n^{3/4}}\cos(4\pi\sqrt{nx}-\frac{\pi}{4}),\\
R_{22}(x)&=\frac{x^{1/4}}{\sqrt 2\pi}\sum_{z <n \leq T}
           \frac{d(n)}{n^{3/4}}\cos(4\pi\sqrt{nx}-\frac{\pi}{4})+O(T^{\varepsilon}).
\end{align*}

Obviously we have
\begin{equation*}
R_{22}(x) \ll x^{1/2+\varepsilon}z^{-1/2}.
\end{equation*}
 For the mean square estimate for $R_{22}(x)$, we have
\begin{equation*}
\int_T^{T+H} R_{22}^2(x)dx \ll T^{1/2+\varepsilon}Hz^{-1/2}+T^{1+\varepsilon}
\end{equation*}
(see e.g. (2.3) of Zhai \cite{Z4}).
Hence
\begin{align}
\int_T^{T+H} R_{22}^4(x)dx
& \ll \max_{T \le x \le T+H}R_{22}^2(x) \int_T^{T+H}R_{22}^2(x)dx \\[1ex]
& \ll T^{1+\varepsilon}z^{-1}(T^{1/2}Hz^{-1/2}+T).\nonumber
\end{align}

For $R_{21}(x)$ we can write
\begin{align*}
R_{21}(x)&=\frac{x^{1/4}}{\sqrt{2}\pi}\sum_{y<n\leq z}\frac{d(n)}{n^{3/4}}\cos(4\pi\sqrt{nx}-\pi/4)\\
&= \sum_{0\leq j\leq J}\frac{x^{1/4}}{\sqrt{2}\pi}\sum_{2^jy<n\leq 2^{j+1}y}
\frac{d(n)}{n^{3/4}}\cos(4\pi\sqrt{nx}-\pi/4).
\end{align*}

By H\"older's inequality we get
\begin{equation*}
|R_{21}(x)|^4\ll TJ^3\sum_{0\leq j\leq J}\left|\sum_{2^jy<n\leq
2^{j+1}y}\frac{d(n)}{n^{3/4}}\cos(4\pi\sqrt{nx}-\pi/4)\right|^4.
\end{equation*}

Correspondingly we have
\begin{align*}
&\int_T^{T+H}|R_{21}(x)|^4dx \\
& \ll TJ^3\sum_{0\leq j\leq J}\int_T^{T+H}\left|\sum_{2^jy<n\leq 2^{j+1}y}
      \frac{d(n)}{n^{3/4}}\cos(4\pi\sqrt{nx}-\pi/4)\right|^4 dx\\
& \ll TJ^4 \int_T^{T+H}\left|\sum_{2^j y< n \leq 2^{j+1}y}\frac{d(n)}{n^{3/4}}
      \cos(4\pi\sqrt{nx}-\pi/4)\right|^4 dx
\end{align*}
for some $0 \leq j \leq J.$ Let $N=2^jy$. % and  $b_n=d(n)(N/n)^{3/4}.$
Thus by Lemma 3.5 we get (note that $J\ll \log T$)
\begin{align}
& \int_T^{T+H}|R_{21}(x)|^4dx
\ll T^{1+\varepsilon}\int_T^{T+H}\left|\sum_{n\sim N} \frac{d(n)}{n^{3/4}} e(2\sqrt{nx})\right|^4dx\\
& \hspace{10mm} =T^{1+\varepsilon} \sum_{n,m,k,l\sim N} \frac{d(n)d(m)d(k)d(l)}{(nmkl)^{3/4}}
                 \int_T^{T+H}e(2\eta x^{1/2})dx \nonumber\\
& \hspace{10mm} \ll \frac{T^{1+\varepsilon}}{N^3} \sum_{n,m,k,l\sim N}
                    \min\left(H,\frac{T^{1/2}}{|\eta|}\right),\nonumber
\end{align}
where we put $\eta=\sqrt n+\sqrt m-\sqrt k-\sqrt l$. By Lemma 3.3 the
contribution of $H$ is (note that in this case $|\eta| \leq \frac{T^{1/2}}{H}$)
\begin{align}
& \ll HT^{1+\varepsilon}N^{-3}(T^{1/2}H^{-1}N^{7/2}+N^2)\\
& \ll T^{3/2+\varepsilon}N^{1/2}+HT^{1+\varepsilon}N^{-1}\nonumber\\
& \ll T^{3/2+\varepsilon}z^{1/2}+HT^{1+\varepsilon}y^{-1}.\nonumber
\end{align}
By Lemma 3.3 again we see that the contribution of $T^{1/2}|\eta|^{-1}$ (in this case
 $|\eta|>T^{1/2}H^{-1}$) is
\begin{align}
& \ll T^{3/2+\varepsilon}N^{-3}\max_{\delta\gg T^{1/2}H^{-1}}
      \frac{1}{\delta}\sum_{\begin{subarray}{c} n,m,k,l \sim N \\ \delta<|\eta|\leq 2\delta
                            \end{subarray}} 1\\
& \ll \max_{\delta\gg T^{1/2}H^{-1}}(T^{3/2+\varepsilon}N^{1/2}
       +T^{3/2+\varepsilon}N^{-1}\delta^{-1})\nonumber\\
& \ll T^{3/2+\varepsilon}N^{1/2}+HT^{1+\varepsilon}N^{-1}\nonumber\\
& \ll T^{3/2+\varepsilon}z^{1/2}+HT^{1+\varepsilon}y^{-1}.\nonumber
\end{align}

Combining (6.16)--(6.19) and noting $z \asymp H^2y^{-2}T^{-1}$ we get
that
\begin{align}
\int_T^{T+H}|R_2(x)|^4dx
& \ll HT^{1+\varepsilon}y^{-1}+T^{3/2+\varepsilon}z^{1/2}\\
& \quad +T^{3/2+\varepsilon}Hz^{-3/2}+T^{2+\varepsilon}z^{-1}\nonumber\\[1ex]
& \ll HT^{1+\varepsilon}y^{-1}+T^{3+\varepsilon}y^3H^{-2}.\nonumber
\end{align}

On the other hand, since $M(x)\ll \frac{U}{T}|R_{1}(x+U)|$ and
$$
\int_T^{T+H}|R_1(x)|^4dx\ll \int_T^{2T}|R_1(x)|^4dx\ll T^2
$$
we have
\begin{align}
\int_T^{T+H}M^4(x)dx & \ll U^4T^{-4}\int_T^{T+H}|R_{1}(x+U)|^4dx \\
                     & \ll U^4T^{-4}T^2 \ll 1 .\nonumber
\end{align}

From (6.20) and (6.21)
\begin{equation} \label{eq622}
\int_T^{T+H}S^4_{2}(x)dx \ll T^{1+\varepsilon}Hy^{-1} +T^{3+\varepsilon}y^3H^{-2}.
\end{equation}
We also remark that by \eqref{jyouken-1} the right hand side of \eqref{eq622} is bounded
above by $T^{\varepsilon}HU^2$.
\bigskip

%***************subsction 6.3**********************************

\subsection{\bf Proof of Theorem 1}\

By the elementary estimate $(a+b)^4-a^4\ll |a^3b|+|b|^4$ we can
write
\begin{equation} \label{delta-4jyou}
(\Delta(x+U)-\Delta(x))^4=S_1^4(x)+O(|S_1^3(x)S_2(x)|+|S_2(x)|^4).
\end{equation}

By (6.15), (6.22)  and  the H\"older's inequality we get
\begin{align} \label{kongouterm}
\int_T^{T+H}|S_1^3(x)S_2(x)|dx
&\ll \left(\int_T^{T+H}S^4_{1}(x)dx\right)^{3/4}\left(\int_T^{T+H}S^4_{2}(x)dx\right)^{1/4} \\[1ex]
&  \ll T^{\varepsilon}(HU^2 )^{3/4}(HTy^{-1})^{1/4}+T^{\varepsilon}(HU^2 )^{3/4}
(T^3y^3H^{-2})^{1/4}\nonumber\\[1ex]
&\ll HT^{3/16+\varepsilon}U^{25/16}+H^{1/4}T^{15/16+\varepsilon}U^{21/16}. \nonumber
\end{align}

Now Theorem 1 follows from (6.10)-(6.14) and (6.22)-(6.24).

Finally we prove Corollary 1. Recall the definition  $\Sigma_6$ in
Section 6.1. Since $U\ll T^{1/2-\varepsilon},$ by Lemma 4.1 we see
that
\begin{equation} \label{eq81}
\Sigma_6\asymp HU^2\log^6 \frac{T}{U^2},
\end{equation}
which combining the formula (2.3) of Theorem 1 gives Corollary 1.

%*****************************  section 7 *************************************************

\section{\bf Proof of Theorem 2}

In this section we prove Theorem 2.   We begin with the well-known
Atkinson's formula (see Ivi\'c \cite[Chapter 15]{I})
\begin{equation}
E(t)=\Sigma_1(t)+\Sigma_2(t)+O(\log^2 t),
\end{equation}
where
\begin{equation}
\begin{split}
& \Sigma_1(t):=\frac{1}{\sqrt 2}\sum_{n\leq N}h(t,n)\cos(f(t,n)),\\
& \Sigma_2(t):=-2\sum_{n\leq N^{\prime}}d(n)n^{-1/2}\left(\log\frac{t}{2\pi n}\right)^{-1}
   \cos\left(t\log\frac{t}{2\pi n}-t+\frac{\pi}{4}\right),\\
& h(t,n):=(-1)^nd(n)n^{-1/2}\left(\frac{t}{2\pi n}+\frac{1}{4}\right)^{-1/4}(g(t,n))^{-1},\\
&g(t,n):={\rm arsinh}((\frac{\pi n}{2t})^{1/2}),\\
&f(t,n):=2tg(t,n)+(2\pi nt+\pi^2n^2)^{1/2}-\pi/4,\\[1ex]
& At\leq N\leq A^{\prime}t, N^{\prime}:=t/2\pi+N/2-(N^2/4+Nt/2\pi)^{1/2},
\end{split}
\end{equation}
where $0<A<A^{\prime}$ are any fixed constants.

We also need the error term $\Delta^{*}(t),$ defined by
\begin{equation}
\Delta^{*}(t): =\frac{1}{2}\sum_{n\leq 4t}(-1)^nd(n)-t(\log t+2\gamma-1).
\end{equation}
$\Delta^{*}(t)$ also has the truncated Voronoi's formula (see Ivi\'c \cite[(15.68)]{I}):
$$
\Delta^{*}(t)=\frac{t^{1/4}}{\sqrt 2\pi}\sum_{n\leq N}\frac{(-1)^n d(n)}{n^{3/4}}
\cos\left(4\pi\sqrt{nt}-\frac{\pi}{4}\right)+O(t^{1/2+\varepsilon}N^{-1/2}),
$$
for $1 \ll N \ll t$.

Let $y:=c(T/U)^{1/4},$ where $c$ is a small positive constant. For $T\ll t\ll T,$ define
\begin{align*}
{\cal E}_1(t): &=\frac{1}{\sqrt 2}\sum_{n\leq y}h(t,n)\cos(f(t,n)),\\
{\cal E}_2(t):&=E(t)-{\cal E}_1(t),\\
{\cal R}_1^{*}(t):&=\frac{t^{1/4}}{\sqrt 2\pi}\sum_{n\leq y}\frac{(-1)^n d(n)}{n^{3/4}}
                    \cos\left(4\pi\sqrt{nt}-\frac{\pi}{4}\right),\\
{\cal R}_2^{*}(t):&=\Delta^{*}(t)-{\cal R}_1^{*}(t).
\end{align*}

{\bf Step 1.} The upper bound of $\int_T^{T+H}{\cal E}_2^4(t)dt$

In this subsection we shall  show that
\begin{equation}
\int_T^{T+H}{\cal E}_2^4(t)dt\ll HT^{1+\varepsilon}y^{-1}.
\end{equation}
Let $z:=(T/y)^{1/2}$ and define
\begin{align*}
{\cal E}_{21}(t): &=\frac{1}{\sqrt 2}\sum_{y<n\leq z}h(t,n)\cos(f(t,n)),\\
{\cal E}_{22}(t): &={\cal E}_{2}(t)-{\cal E}_{21}(t).
\end{align*}

Let $1/4<\theta< 1/3$ be a constant such that $E(t)\ll t^\theta.$
Following Ivi\'c \cite{I}, the second author proved in \cite{Z4} that
\begin{equation*}
\int_{T}^{T+H}|\Delta(t)|^Adt\ll HT^{A/4+\varepsilon}
\end{equation*}
holds for $H \geq T^{1+\theta(A-2)-A/4}$ if
$2<A<2\theta/(\theta-1/4).$  By the same argument we can show that the estimate
\begin{equation*}
\int_{T}^{T+H}|E(t)|^Adt\ll HT^{A/4+\varepsilon}
\end{equation*}
holds for $H\geq T^{1+\theta(A-2)-A/4}$ if
$2<A<2\theta/(\theta-1/4).$ A slight modification  shows that if
$H\geq T^{1+\theta(A-2)-A/4}$ and $2<A<2\theta/(\theta-1/4),$ then
\begin{equation}
\int_{T}^{T+H}|{\cal E}_{22}(t)|^Adt\ll HT^{A/4+\varepsilon}.
\end{equation}
 Similar to  the formula (2.8) of Zhai \cite{Z4},  we can
easily show that
\begin{equation}
\int_{T}^{T+H}|{\cal E}_{22}(t)|^2dt\ll HT^{1/2}z^{-1/2}\log^3 T.
\end{equation}
We omit the proofs of the above formulas. From (7.5), (7.6) and
H\"older's inequality we get
\begin{equation*}
\int_{T}^{T+H}|{\cal E}_{22}(t)|^4dt\ll
HT^{1+\varepsilon}z^{-\frac{A-4}{2(A-2)}}
\end{equation*}
holds for  $H\geq T^{1+\theta(A-2)-A/4}$ if
$4<A<2\theta/(\theta-1/4).$ Taking $\theta=72/227$ and $A=8$ we get
that the estimate
\begin{align}
\int_{T}^{T+H}|{\cal E}_{22}(t)|^4dt
& \ll HT^{1+\varepsilon}z^{-1/3}\ll HT^{1+\varepsilon}(T/y)^{-1/6}\\
& \ll HT^{1+\varepsilon}y^{-1}\nonumber
\end{align}
holds in the range $T^{205/227}\ll H\leq T$ if noting that
$T^{3/7}\ll U\ll T^{1/2}.$

Now we estimate the integral $\int_T^{T+H}|{\cal E}_{21}(t)|^4dt.$
Similar to (6.17) we have
\begin{align}
\log^{-4} T&\int_T^{2T}|{\cal E}_{21}(t)|^4dt
\ll \int_T^{2T}\left|\sum_{n\sim N}h(t,n)e(f(t,n))\right|^4dt\\
& =\sum_{N<n,m,k,l\leq
2N}\int_T^{2T}H(t;n,m,k,l)e(F(t;n,m,k,l))dt\nonumber
\end{align}
holds for some $y<N<z$, where
\begin{align*}
H(t;n,m,k,l)&=h(t,n)h(t,m)h(t,k)h(t,l),\\
F(t;n,m,k,l)&=f(t,n)+f(t,m)-f(t,k)-f(t,l).
\end{align*}
From (7.2) it is easy to check that
\begin{equation}
\begin{split}
h(t,n)&=\frac{2^{3/4}}{\pi^{1/4}}\frac{(-1)^nd(n)}{n^{3/4}}t^{1/4}(1+O(n/t)),\\
f(t,n)&=2^{3/2}(\pi nt)^{1/2}-\pi/4+O(n^{3/2}t^{-1/2}),\\
f^{\prime}(t,n)&=2^{1/2}(\pi n)^{1/2}t^{-1/2}+O(n^{3/2}t^{-3/2}).
\end{split}
\end{equation}
Let $\eta=\sqrt{n}+\sqrt{m}-\sqrt{k}-\sqrt{l}.$ From the third formula of (7.9) we get
$$
F^{\prime}(t;n,m,k,l)=(2\pi)^{1/2}\eta t^{-1/2}+O(z^{3/2}t^{-3/2}).
$$
Let $C>0$ be a large real constant such that if $|\eta|\geq Cz^{3/2}T^{-1}$
then  $|F^{\prime}(t;n,m,k,l)|\gg |\eta|T^{-1/2}.$

If $|\eta|\leq Cz^{3/2}T^{-1},$ then by  the trivial estimate and
Lemma 3.3 we get
\begin{align}
&\int_T^{T+H}\sum_{\begin{subarray}{c} N<n,m,k,l\leq 2N \\ |\eta|\leq Cz^{3/2}T^{-1}\end{subarray}}
   H(t;n,m,k,l)e(F(t;n,m,k,l))dt\\
& \ll \frac{HT^{1+\varepsilon}}{N^3}
\sum_{\begin{subarray}{c} N<n,m,k,l\leq 2N \\ |\eta|\leq Cz^{3/2}T^{-1} \end{subarray}} 1 \nonumber  \\
& \ll \frac{HT^{1+\varepsilon}}{N^3}(Cz^{3/2}T^{-1}N^{7/2}+N^2)N^\varepsilon \nonumber\\
& \ll HT^{1+\varepsilon}y^{-1}+Hz^{2}T^{\varepsilon} \nonumber \\
& \ll HT^{1+\varepsilon}y^{-1}.\nonumber
\end{align}

Now suppose $|\eta|>Cz^{3/2}T^{-1}.$  By Lemma 3.5 and 3.3  we get
\begin{eqnarray}
&&\int_T^{T+H}\sum_{\stackrel{N<n,m,k,l\leq 2N}{|\eta|>
Cz^{3/2}Y^{-1}}}H(t;n,m,k,l)e(F(t;n,m,k,l))dt\\
&&\ll \frac{T^{3/2+\varepsilon}}{N^3}\sum_{\stackrel{N<n,m,k,l\leq
2N}{|\eta|> Cz^{3/2}T^{-1}}}\frac{1}{|\eta|} \ll
\frac{T^{3/2+\varepsilon}}{N^3\delta}\sum_{\stackrel{N<n,m,k,l\leq
2N}{|\eta|\sim \delta
\gg z^{3/2}T^{-1}}}1\nonumber\\
&&\ll   \frac{T^{3/2+\varepsilon}}{N^3\delta} (\delta
N^{7/2}+N^2)N^\varepsilon\nonumber\ll
\frac{HT^{1+\varepsilon}}{y}+\frac{T^{3/2+\varepsilon}}{N\delta}\nonumber\\
&&\ll
\frac{HT^{1+\varepsilon}}{y}+\frac{T^{5/2+\varepsilon}}{Nz^{3/2}}\ll
\frac{HT^{1+\varepsilon}}{y}+\frac{T^{7/4+\varepsilon}}{y^{1/4}}\nonumber\\&&
\ll \frac{HT^{1+\varepsilon}}{y} .\nonumber
\end{eqnarray}

From (7.8) , (7.10) and (7.11) we get
\begin{eqnarray}
\int_T^{T+H}|{\cal E}_{21}(t)|^4dt\ll
HT^{1+\varepsilon}y^{-1}+T^{3/2+\varepsilon}z^{1/2},
\end{eqnarray}
which combining (7.7) gives (7.4).

\bigskip

{\bf Step 2.} The upper bound of $\int_T^{T+H}\left({\cal
E}_1(t)-2\pi{\cal R}_1^{*}(t/2\pi)\right)^{4}dt$

It is easy to see that
\begin{eqnarray*}
2\pi R_1^{*}(\frac{t}{2\pi})=\frac{(2t)^{1/4}}{\pi^{1/4}}\sum_{n\leq
y}\frac{(-1)^n d(n)}{n^{3/4}}\cos\left(2^{3/2}(\pi
nt)^{1/2}-\frac{\pi}{4}\right).
\end{eqnarray*}
Thus we can  write
\begin{eqnarray*}
&&{\cal E}_1(t)-2\pi R_1^{*}(\frac{t}{2\pi})=S_3(t)+S_4(t),\\
&&S_3(t)=\sum_{n\leq y}h_1(t,n)\cos(f(t,n)),\\
&&h_1(t,n)=\frac{h(t,n)}{\sqrt
2}-\frac{(2t)^{1/4}}{\pi^{1/4}}\frac{(-1)^n
d(n)}{n^{3/4}},\\
&&S_4(t)=\sum_{n\leq y}\frac{(2t)^{1/4}}{\pi^{1/4}}\frac{(-1)^n
d(n)}{n^{3/4}}\left(\cos(f(t,n))-\cos\left(2^{3/2}(\pi
nt)^{1/2}-\frac{\pi}{4}\right)\right).
\end{eqnarray*}
From (7.9) we have $h_1(t,n)\ll d(n)n^{1/4}t^{-3/4}$ and $S_3(t)\ll
y^{5/4}t^{-3/4}\log y.$ Thus
\begin{equation}
\int_T^{T+H}|S_3(t)|^4dt\ll y^5HT^{-3}\log^4 y\ll 1.
\end{equation}

By the simple relation
\begin{equation}
\cos(u-v)-\cos(u+v)=2\sin u\sin v
\end{equation}
we can write
\begin{eqnarray*}
&&S_4(t)=\sum_{n\leq y}h_2(t,n)\sin(f_1(t,n)),\\
&&h_2(t,n)=\frac{(2t)^{1/4}}{\pi^{1/4}}\frac{(-1)^n d(n)}{n^{3/4}}
\sin \frac{f(t,n)-2^{3/2}(\pi nt)^{1/4}+\pi/4}{2},\\
&&f_1(t,n)=\frac{f(t,n)+2^{3/2}(\pi nt)^{1/4}-\pi/4}{2}.
\end{eqnarray*}

It is easy to check that
$$\frac{f(t,n)-2^{3/2}(\pi
nt)^{1/4}+\pi/4}{2}=\beta_3n^{3/2}t^{-1/2}+\beta_5n^{5/2}t^{-3/2}+\cdots,$$
where $\beta_3,\beta_5,\cdots$ are real constants. So for each
$n\leq y,$ $h_2(t,n)$ is a monotonic function of $t$ and
$h_2(t,n)\ll d(n)n^{3/4}t^{-1/4}.$ By (7.14) again we can write
\begin{eqnarray*}
|S_4(t)|^2&&=\sum_{n,m\leq
y}h_2(t,n)h_2(t,m)\sin(f_1(t,n))\sin(f_2(t,m))\\
&&=S_5(t)+S_6(t)+S_7(t),\\
S_5(t)&&=\sum_{n\leq y}h_2^2(t,n)\sin^2(f_1(t,n)),\\
S_6(t)&&=\frac{1}{2}\sum_{\stackrel{n,m\leq y}{n\not= m}}
h_2(t,n)h_2(t,m)\cos(f_1(t,n)-f_1(t,m)),\\
S_7(t)&&=\frac{1}{2}\sum_{\stackrel{n,m\leq y}{n\not= m}}
h_2(t,n)h_2(t,m)\cos(f_1(t,n)+f_1(t,m)).
\end{eqnarray*}
Trivially $S_5(t)\ll y^{5/2}t^{1/2}\log^3 y,$ which implies
\begin{equation}
\int_T^{T+H}S_5(t)dt\ll y^{5/2}HT^{-1/2}\log^3 y.
\end{equation}
From the third formula of (7.9) it is easy to check that
$$|f_1^{\prime}(t,n)-f_1^{\prime}(t,m)|\gg |\sqrt n-\sqrt m|T^{-1/2}.$$
By Lemma 3.5 we have
\begin{eqnarray}
&&\int_T^{T+H}S_6(t)dt\ll \sum_{n\not=
m}\frac{d(n)d(m)(nm)^{3/4}}{|\sqrt n-\sqrt m|}\\
&&\hspace{5mm}\ll \sum_{|\sqrt n-\sqrt m|\geq
(nm)^{1/4}/100}\frac{d(n)d(m)(nm)^{3/4}}{|\sqrt n-\sqrt
m|}\nonumber\\
&&\hspace{8mm} +\sum_{|\sqrt n-\sqrt m|<
(nm)^{1/4}/100}\frac{d(n)d(m)(nm)^{3/4}}{|\sqrt n-\sqrt
m|}\nonumber\\
&&\hspace{5mm}\ll \sum_{n,m\leq y}d(n)d(m)(nm)^{1/2}+ \sum_{n\asymp
m\leq y}\frac{d(n)d(m)(nm)^{3/4}}{n^{-1/2}| n-m|}\nonumber\\
&&\hspace{5mm}\ll y^3\log^2 y+\sum_{n\asymp m\leq
y}\frac{d^2(n)n^2}{| n-m|}\ll y^3\log^4 y.\nonumber
\end{eqnarray}
Similarly
\begin{equation}
\int_T^{T+H}S_7(t)dt\ll \sum_{n\not=
m}\frac{d(n)d(m)(nm)^{3/4}}{|\sqrt n+\sqrt m|}\ll y^3\log^2 y.
\end{equation}

From (7.15)-(7.17) we get
\begin{eqnarray*}
\int_T^{T+H}|S_4(t)|^2dt\ll y^{5/2}HT^{-1/2}\log^3 y+y^3\log^4 y\ll
y^{5/2}HT^{-1/2}\log^3 y,
\end{eqnarray*}
which combining the trivial estimate $S_4(t)\ll y^{7/4}t^{-1/4}\log
y$ gives
\begin{eqnarray}
\int_T^{T+H}|S_4(t)|^4dt\ll y^{6}HT^{-1}\log^3 y\ll HTy^{-1}.
\end{eqnarray}
So from (7.13) and (7.18) we get
\begin{eqnarray}
&&\int_T^{T+H}|{\cal E}_1(t)-2\pi{\cal R}_1^{*}(t/2\pi)|^4dt \ll
HTy^{-1}.
\end{eqnarray}

{\bf Step 3.} Proof of Theorem 2

For $T\leq t\leq T+H,$ we can write
\begin{eqnarray}
&&\hspace{6mm}E(t+U)-E(t)\\&&={\cal E}_1(T+U)-{\cal E}_1(t)+{\cal
E}_2(T+U)-{\cal
E}_2(t)\nonumber\\
&&=2\pi{\cal R}_1^{*}(\frac{t+U}{2\pi})-2\pi{\cal
R}_1^{*}(\frac{t}{2\pi})+{\cal E}_1(T+U)-2\pi{\cal
R}_1^{*}(\frac{t+U}{2\pi})\nonumber\\
&&\hspace{10mm}-{\cal E}_1(t)+2\pi{\cal
R}_1^{*}(\frac{t}{2\pi})+{\cal
E}_2(T+U)-{\cal E}_2(t)\nonumber\\
&&=S_8(t)+S_9(t)\nonumber
\end{eqnarray}
say, where
\begin{eqnarray*}
S_8(t):&&=2^{1/4}\pi^{-1/4}t^{1/4}\sum_{n\leq y}\frac{(-1)^n
d(n)}{n^{3/4}}\cos\left(2^{3/2}\pi^{1/2}\sqrt{nt}-\frac{\pi}{4}\right),\\
S_9(t):&&={\cal E}_1(T+U)-2\pi{\cal R}_1^{*}(\frac{t+U}{2\pi})
-{\cal E}_1(t)+2\pi{\cal
R}_1^{*}(\frac{t}{2\pi})\\&&\hspace{15mm}+{\cal E}_2(T+U)-{\cal
E}_2(t)+M^{*}(t),\\
M^{*}(t):&&=2^{1/4}\pi^{-1/4}\left((t+U)^{1/4}-t^{1/4}\right)\\
&&\hspace{5mm}\times \sum_{n\leq y}\frac{(-1)^n
d(n)}{n^{3/4}}\cos\left(2^{3/2}\pi^{1/2}\sqrt{n(t+U)}-\frac{\pi}{4}\right)\\
&&\ll UT^{-1}|{\cal R}_1^{*}(t+U)|.
\end{eqnarray*}
Similar to (6.21), we have
\begin{eqnarray*}
\int_T^{T+H}M^{*4}(t)dt\ll 1,
\end{eqnarray*}
which combining (7.4) and (7.19) gives
\begin{equation}
\int_T^{T+H}S_9^4(t)dt\ll HT^{1+\varepsilon}y^{-1}.
\end{equation}
For $S_8(t),$ similar to (6.10) under the condition (2.1) we have
the asymptotic formula
\begin{eqnarray}
&&\int_{T}^{T+H}S_{8}^4(t)dt=\frac{12}{\pi}\sum_{\stackrel{\sqrt{n_1}+\sqrt{n_2}
=\sqrt{n_3}+\sqrt{n_4}}{n_j\leq
y}}\frac{d(n_1)d(n_2)d(n_3)d(n_4)}{(n_1n_2n_3n_4)^{3/4}}\\
&&\hspace{10mm}\times \int_T^{T+H}t\prod_{j=1}^4 \sin \frac{
U\sqrt{\pi n_j}}{(2t)^{1/2}}dx+O(HU^2T^{-1/2}\log^4 T)\nonumber\\
&&\hspace{15mm}+O(HT^{3/16+\varepsilon}U^{25/16}).\nonumber
\end{eqnarray}
  Now  Theorem 2 follows from (7.20)--(7.22).

\begin{flushleft}
Y. Tanigawa\\

Graduate School of Mathematics, Nagoya University, \\
Chikusa-ku, Nagoya 464-8602, Japan\\
E-mail: tanigawa@math.nagoya-u.ac.jp\\

\bigskip

Wenguang Zhai \\
School of Mathematical Sciences, Shandong Normal University, \\
Jinan 250014, P. R. China\\
E-mail: zhaiwg@hotmail.com
\end{flushleft}

\end{document}